\documentclass{article}

\usepackage{graphicx}
   \newtheorem{Tw}{Theorem}[section]
   
   \newtheorem{De}[Tw]{Definition}
   \newtheorem{Le}[Tw]{Lemma}

 \let\pa=\partial
   \newcommand{\qed}{$\quad\Box$\par\medskip}

   \newcommand{\init}{\mathrm{in}}

   \newcommand{\supp}{\mathrm{supp}}
   \newcommand{\convex}{\mathrm{convex}}
   
   \newcommand{\Pole}{\mathrm{Area}}

   \newcommand{\grad}{\mathrm{grad}\,}
   
   \newcommand{\Rn}{\mathbf{R}}
   \newcommand{\Cn}{\mathbf{C}}
   \newcommand{\Nn}{\mathbf{N}}

\usepackage{amssymb}

\def\e{\kern+.5ex\lower.42ex\hbox{$\scriptstyle \iota$}\kern-1.10ex e}

\title{\large NEARLY IRREDUCIBILITY OF POLYNOMIALS AND THE NEWTON DIAGRAMS
\footnotetext{
\phantom{ssssi}Key words and phrases: irreducibility of polynomials, Newton diagram, \\
\phantom{ssssi}Newton polygon, plane algebraic curve.
}
}

\begin{document}
\author{\large Mateusz Masternak}
\date{}
\maketitle

\begin{abstract}
Let $f$ be a polynomial in two complex variables. We say that $f$ is nearly irreducible if any two nonconstant polynomial factors of $f$  have a common zero in $\Cn^2$. In the paper we give a~criterion of nearly irreducibility for a~given polynomial $f$ in terms of its Newton diagram.

\end{abstract}
\vspace{3ex}

\section{Introduction}
Let $f(X,Y)=\sum c_{\alpha\beta}X^{\alpha}Y^{\beta}\in\Cn[X,Y]$ be a nonzero
  polynomial of positive degree. We say that the polynomial $f$ is
  \emph{quasi--convenient} if
  $c_{\alpha0}\neq0$ and $c_{0\beta}\neq0$ for some integers
  $\alpha,\beta\geq0$.  Let $\supp f:=\{(\alpha,\beta)\in\Nn^2:\;c_{\alpha\beta}\neq 0\}$.
  We define
$$
\Delta_{\infty}(f) := \convex(\{(0,0)\}\cup\supp f).
$$
  The polygon  $\Delta_{\infty}(f)$ is called  \emph{the Newton  diagram at infinity} of the polynomial~$f$.
  
  For any nonzero  vector $\vec w=[p,q]$ of   the real plane~$\Rn^2$ we put
  $$ \init(f,\vec w)(X,Y):= \sum_{p\alpha+q\beta=d_{\vec w}(f)}
   c_{\alpha\beta}X^{\alpha}Y^{\beta},
$$
  where 
  $$ d_{\vec w}(f)=\max{\{p\alpha+q\beta: (\alpha,\beta)\in\supp f\}} . 
  $$
   
 We call a quasi--convenient polynomial $f(X,Y)\in\Cn[X,Y]$  \emph{nondegenerate at infinity} if for every real vector $\vec w=[p,q]$  such that $p>0$   or $q>0$ the system of equations 
 $$ \init(f,\vec
  w)(X,Y) = \frac{\pa}{\pa X}\init(f,\vec w)(X,Y) = \frac{\pa}{\pa
  Y}\init(f,\vec w)(X,Y) = 0 
  $$has no solutions in 
  $\Cn^{*}\times\Cn^{*}$, where $\Cn^*=\Cn\setminus\{0\}$. \newpage
  
\begin{De}\label{De:1.1}\sloppy
  A polynomial $f(X,Y)\in\Cn[X,Y]$ of a positive degree is nearly irreducible
  if any two nonconstant polynomial factors
  $g(X,Y)$, \mbox{$h(X,Y)\in\Cn[X,Y]$~of~$f(X,Y)$} have a common zero
  in $\Cn^2$.
\end{De}

Note that every nearly irreducible polynomial has a connected
  zero-set. Note that nearly irreducible polynomial  may be reducible (e.g.~$f=XY$).
  It is easy to check that if $f$ is nearly irreducible and
  $\grad f =(\frac{\pa f}{\pa X}, \frac{\pa f}{\pa Y})\neq0$
   on the curve $f(X,Y)=0$ then $f$ is
  irreducible~(see \cite{RST}).

The notion of nearly irreducibility of polynomials in two variables was
  introduced in~\cite{AR} by S.~Abhyankar and L.~A.~Rubel in
  connection with research of these authors on irreducibility of
  polynomials of the form $f(X)-g(Y)$. The main result of~\cite{AR} was reproved by L.~A.~Rubel, A.~Shinzel and
  H.~Tverberg in~\cite{RST}. Afterwards A.~P\l{}oski 
  generalized the result of Abhyankar and Rubel by using  the Newton diagram of a given polynomial (see  ~\cite{P2}, Theorem 2) which is  Theorem \ref{Tw:1.1} in this note. 
  \begin{Tw}[\cite{P2}, Theorem 2]\label{Tw:1.1}\
Let $f$ be a quasi-convenient polynomial such that
\begin{itemize}
\item[(1)] $f$ is nondegenerate at infinity,
\item[(2)] every face of the polygon $\Delta_\infty(f)$ not included in  coordinate axes has a~negative slope, (i.e. it is a segment included in the straight line of the form $p\alpha+q\beta=r$ for some $p,q>0$).
\end{itemize}
Then the polynomial $f$ is nearly irreducible.
\end{Tw}

  Our theorem (Theorem \ref{Tw:1.2}) generalizes result of P\l{}oski. We state

\begin{Tw}\label{Tw:1.2}\
Let $f$ be a quasi-convenient polynomial such that
\begin{itemize}
\item[(1)] $f$ is nondegenerate at infinity,
\item[(2)] if $\vec{w}=[p,q]$ is a nonzero vector such that $pq\leqslant 0$ then the system of equations
$\init(f,\vec{w})(X,Y)\!=\!\init(f,-\vec{w})(X,Y)\!=\!0$ has no solutions in $\Cn^*\!\times\!\Cn^*$.
\end{itemize}
Then the polynomial $f$ is nearly irreducible.
\end{Tw}

 In comparison to Theorem  \ref{Tw:1.1} in Theorem \ref{Tw:1.2} there is no restrictions on the shape of the polygon $\Delta_\infty (f)$. 
The proof of Theorem \ref{Tw:1.2}, based on the Kouchnirenko-Bernstein Theorem is given in Section \ref{s:3}.\medskip

\noindent{\bf Remark.} If there is no a pair of parallel faces of the polygon $\Delta_\infty(f)$ then for any $\vec{w}\neq\vec{0}$ at least one of the polynomials $\init(f,\vec{w})(X,Y)$ or $\init(f,-\vec{w})(X,Y)$ is a~monomial  and then the condition (2) in our theorem trivially holds, so Theorem \ref{Tw:1.2} implies Theorem \ref{Tw:1.1} (see Example 3).

 The examples presented below show that the assumption (2) in Theorem \ref{Tw:1.2}  is essential.\newpage

\medskip
\textbf{Example~1.}
  Let us consider the polynomial
$$ f(X,Y)=X^3Y^2+X^2Y^3-X-Y=(XY-1)(XY+1)(X+Y).
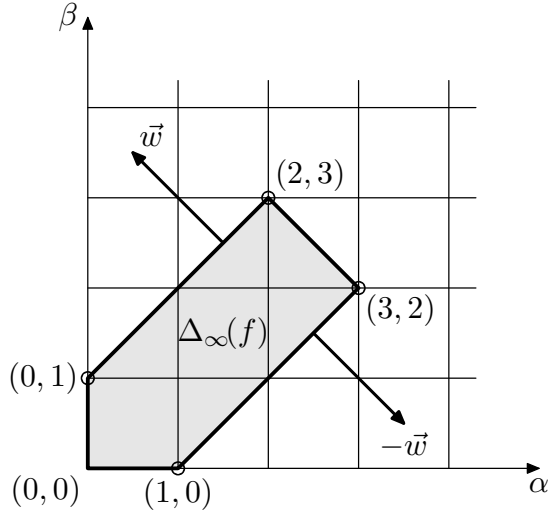
$$

\begin{center}\noindent\includegraphics[scale=1.2]{rys.1}
\end{center}

It is easily seen that the polynomial $f$ is nondegenerate
     at infinity and that it is not nearly irreducible.
    Note that the condition (2) of Theorem \ref{Tw:1.2}  is not satisfied. Namely if $\vec{w}=[-1,1]$ then $\init(f,\vec{w})(X,Y)=-Y(1+X^2Y^2)$ and $\init(f,-\vec{w})(X,Y)=-X(1+X^2Y^2)$.
  The system $\init(f,\vec{w})(X,Y)=\init(f,-\vec{w})(X,Y)=0$ has a~solution in
  $\Cn^{*}\times\Cn^{*}$.
\medskip
\medskip
\medskip

\textbf{Example~2.}  Let $f(X,Y)=(X-1)(X+1)(X+Y)=X^3+X^2Y-X-Y$.
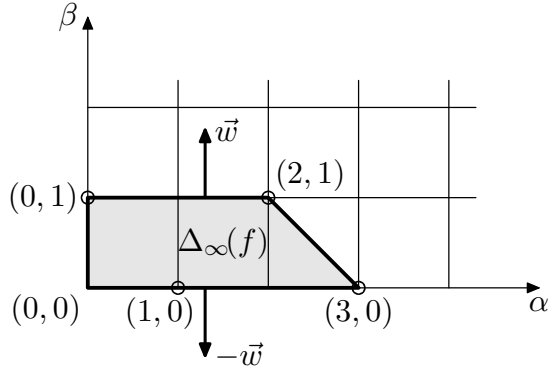
\medskip

\begin{center}\noindent\includegraphics[scale=1.2]{rys.2}
\end{center}

\sloppy The polynomial $f$ is nondegenerate at infinity and obviously $f$ is not nearly irreducible.
 The assumption (2)  of Theorem \ref{Tw:1.2}  does not  hold because if $\vec{w}=[0,1]$  then $\init(f,\vec{w})(X,Y)=Y(X^2-1)$ and  $\init(f,-\vec{w})(X,Y)=X(X^2-1)$  have a common zero
  in   $\Cn^{*}\times\Cn^{*}$.
  Note that for any $c\neq 0$ the polynomial $f(X,Y)+c$ satisfies (2), so it is nearly irreducible.\medskip
 \newpage
  
  \textbf{Example 3.}  Let $f(X,Y)=Y+X^2+XY^3+X^3Y^4+X^5Y^3$.
  
  
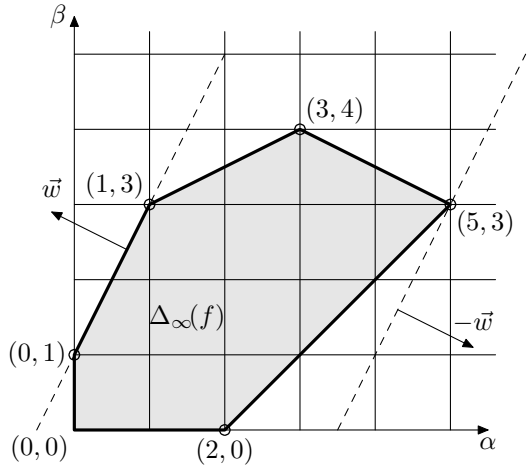
\begin{center}\includegraphics[scale=1]{rys.3}
\end{center}

The polynomial $f$ is nondegenerate at infinity and $\init(f,\vec{w})(X,Y)$  or $\init(f,-\vec{w})(X,Y)$ is a monomial for any $\vec{w}\neq\vec{0}$, hence the polynomial $f$ is nearly irreducible.

\section{Kouchnirenko-Bernstein Theorem}\label{s:2}

  The famous B\'{e}zout theorem for affine curves states that two polynomia
  equations of given degree $m$, $n>0$
  have at most $mn$ common solutions provided that their number is finite.
  If additionally their Newton
  diagrams at infinity are known then we can give more precise estimation. Namely, we may replace
  the product $mn$ by the Minkowski's mixed area of these diagrams.
  Such results were proved in Kouchnirenko and Bernstein's papers
  in 1970s~\cite{Ku1,Ku3,Ku4,B}. See also \cite{A,BKuKh, Kh1, Kh2}. Focusing only on two-dimensional case much more precise results are possible.\\
 
Let $f(X,Y)$, $g(X,Y)\in\Cn[X,Y]$ be polynomials  of positive degrees.
 If $P=(a,b)\in\Cn^2$ is a solution of the system
$$  f(X,Y)=0, \quad
    g(X,Y)=0
$$
 then the symbol $(f,g)_P$\glossary{$(f,g)_P$} denotes the
intersection multiplicity\index{krotno\s{}\c{}}. We use the definition of the intersection multiplicity
 as in~\cite{F}. We have $(f,g)_P<{+}\infty$ if and only if
 $P$ is an isolated solution of the given system.

The pair of quasi--convenient polynomials $(f,g)$ is
           \emph{nondegenerate at infinity}
           if for any real vector
           $\vec w=[p,q]$ such that $p>0$ or
           $q>0$ the system of equations
           $\init(f,\vec w)(X,Y)=\init(g,\vec w)(X,Y)=0$  has no
           solutions in $\Cn^{*}\times\Cn^{*}$.

 For a pair of quasi--convenient polynomials $f(X,Y)$,
  $g(X,Y)\in\Cn[X,Y]$  we denote
$ \nu_{\infty}(f,g) :=         \glossary{$\nu_{\infty}(f,g)$}
   \Pole\Delta_{\infty}(fg)-\Pole\Delta_{\infty}(f)-\Pole\Delta_{\infty}(g)$.\medskip
   
   Let us present a useful version of the  Kouchnirenko--Bernstein Theorem in two-dimensional case.

\begin{Tw}[Kouchnirenko--Bernstein]\label{Tw:2.1}\sloppy
  Let  polynomials  $f(X,Y)$, 
  \mbox{$g(X,Y)\in\Cn[X,Y]$} be quasi--convenient. It holds
\begin{itemize}
\item[(1)] if  $f$ and $g$ are coprime then
           $\quad \sum_{P\in\Cn^2}(f,g)_P \leq \nu_{\infty}(f,g) $,
\item[(2)] $\sum_{P\in\Cn^2}(f,g)_P=\nu_{\infty}(f,g)$  if and only if 
            the pair $(f,g)$ is nondegenerate at infinity.
\end{itemize}
\end{Tw}

  The first proof of this theorem (in multi-dimensional case) was given  by Kouchnirenko in~\cite{Ku1}
  with additional assumption that the polynomials $f$ and $g$
  have  identical Newton diagrams at infinity.

  The original  Bernstein Theorem was formulated for Laurent polynomials without mentioned Kouchnirenko's assumption. Theorem \ref{Tw:1.2} follows from its local version due to Kouchnirenko (i.e. estimation of the intersection multiplicity  of plane curves given in terms of their local Newton diagrams, see~\cite{Ku1,AJ,P1, P3, Ma2}) and from B\'{e}zout Theorem for projectives curves.

\section{Proof of Theorem \ref{Tw:1.2}}\label{s:3}

  The proof of our theorem needs two lemmas.

\begin{Le}\label{Le:3.1}
  Let $f(X,Y)$ be a quasi--convenient polynomial nondegenerate at
  infinity and let $g(X,Y)$ and $h(X,Y)$ be two coprime divisors
  in $\Cn[X,Y]$.  Then the pair $(g,h)$ is nondegenerate at infinity.
\end{Le}
\textbf{Proof of Lemma \ref{Le:3.1}.}\\
  Since $f(X,0)f(0,Y)\neq0$ thus $g(X,0)g(0,Y)h(X,0)h(0,Y)\neq0$ in $\Cn[X,Y]$.
  Therefore the polynomials $g$ and $h$ are quasi--convenient.
  Let us suppose, contrary to our claim, that the pair $(g,h)$ is 
  degenerate at infinity.  By definition 
  there exists a real vector $\vec
  w=[p,q]$, where $p>0$ or $q>0$, such that $\init(g,\vec
  w)(x,y)=\init(h,\vec w)(x,y)=0$ for some $(x,y)\in\Cn^{*}\times\Cn^{*}$.
  Since $g(X,Y)$ and $h(X,Y)$ are coprime divisors of the polynomial
  $f(X,Y)$ then there exists a~polynomial $P(X,Y)$ such that
  $f(X,Y)=g(X,Y)h(X,Y)P(X,Y)$. Let us note that
  $$ 
  \init(f,\vec w)(X,Y) = \init(g,\vec w)(X,Y)\init(h,\vec
   w)(X,Y)\init(P,\vec w)(X,Y), 
   $$ 
 hence
    $$
\init(f,\vec w)(x,y) = \frac{\pa}{\pa X}\init(f,\vec w)(x,y) = \frac{\pa}{\pa
   Y}\init(f,\vec w)(x,y) = 0. 
   $$ 
    The above equalities contradict
  nondegeneracy at infinity of the polynomial $f$.
   \qed
\newpage
  
\begin{Le}\label{Le:3.2}
  If the polynomials $f$, $g$  $\in\,\Cn[X,Y]$   of positive degrees are quasi--convenient
  then
\begin{itemize}
\item[(1)] $\nu_{\infty}(f,g) \geq 0$,
\item[(2)] $\nu_{\infty}(f,g)=0$ if and only if  the diagrams
           $\Delta_{\infty}(f)$ and $\Delta_{\infty}(g)$ form segments
           included in the same straight line passing through the origin.
\end{itemize}
\end{Le}

  In the proof of   Lemma \ref{Le:3.2} we need the following
  Brunno--Minkowski inequality (see \cite{We}, Theorem~6.5.3):

\begin{Tw}\label{Tw:3.3}
  If $A$ and $B$ are nonempty and measurable subsets of $\Rn^2$
   then
$$ (\Pole(A+B))^{1/2} \geq (\Pole A)^{1/2} + (\Pole B)^{1/2}.
$$
\end{Tw}
\textbf{Proof of Lemma~\ref{Le:3.2}.}
  Note that
  $\Delta_{\infty}(fg)=\Delta_{\infty}(f)+\Delta_{\infty}(g)$. Using  Brunno--Minkowsky inequality for the sets
  $A=\Delta_{\infty}(f)$ and $B=\Delta_{\infty}(g)$ we have
$$ (\Pole\Delta_{\infty}(fg))^{1/2} \geq
   (\Pole\Delta_{\infty}(f))^{1/2} + (\Pole\Delta_{\infty}(g))^{1/2},
$$
  hence
$$ \Pole\Delta_{\infty}(fg) \geq \Pole\Delta_{\infty}(f) +
   \Pole\Delta_{\infty}(g) +
   2[(\Pole\Delta_{\infty}(f))(\Pole\Delta_{\infty}(g))]^{1/2}.
$$
  This proves~(1).

  Suppose now that in~(1) the equality holds. Last inequality implies
  that  $\Pole\Delta_{\infty}(f)=0$ or
  $\Pole\Delta_{\infty}(g)=0$.  Suppose, without  loss of generality,  that 
  ~$\Pole\Delta_{\infty}(f)=0$.  Since the set $\Delta_{\infty}(f)$ is convex,   $(0,0)\in\Delta_{\infty}(f)$ and $\deg
  f>0$ we get that  $\Delta_{\infty}(f)$ is a segment included in a straight line
  passing through the origin.  Moreover $$
\Pole\Delta_{\infty}(fg) = \Pole\Delta_{\infty}(g).  $$
 It is easy to check  that  the diagram $\Delta_{\infty}(g)$ does not contain a point different from the
  origin not belonging to the straight line including 
  $\Delta_{\infty}(f)$. Otherwise we would have
  $\Pole\Delta_{\infty}(fg)>\Pole\Delta_{\infty}(g)$. The last observation
  proves~(2).\qed
  
\noindent\textbf{Proof of Theorem~\ref{Tw:1.2}.}
  Let us suppose, contrary to our claim, that there exist polynomials $g(X,Y)$, $h(X,Y)\in\Cn[X,Y]$ of
  positive degrees  being divisors of the
  polynomial $f(X,Y)$  such that $$ \sum_{P\in\Cn^2}(g,h)_P=0.
  $$ Obviously the polynomials $g(X,Y)$ and $h(X,Y)$ are coprime and they are
  quasi--convenient. From 
  Lemma~\ref{Le:3.1} it follows that the pair~$(g,h)$ is nondegenerate
  at infinity. Using now
  Kouchnirenko--Bernstein Theorem (Theorem~\ref{Tw:2.1}) we state
  that $$ \nu_{\infty}(g,h)=0.  $$
  
   Therefore,
  Lemma~\ref{Le:3.2} implies that the diagrams $\Delta_{\infty}(g)$ and
  $\Delta_{\infty}(h)$ are segments included at the same straight line $p\alpha+q\beta=0$, where $\vec{w}=[p,q]\neq\vec{0}$ and $pq\leqslant 0$.  So we have 
 $$
 \begin{array}{l}
  \init(g,\vec w)(X,Y)= \init(g,-\vec w)(X,Y) = g(X,Y),\\
  \init(h,\vec w)(X,Y) =  \init(h,-\vec w) (X,Y)= h(X,Y).  
  \end{array}
$$ 
There exists a polynomial $P(X,Y)$ such that $$
  f(X,Y)=g(X,Y)h(X,Y)P(X,Y), $$
  hence 
$$
\begin{array}{l}
 \init(f,\vec w)(X,Y)= \init(g,\vec w)(X,Y)\init(h,\vec w)(X,Y)  \init(P,\vec w)(X,Y),\\
 \init(f,{-}\vec w)(X,Y) = \init(g,-\vec w)(X,Y)\init(h,-\vec w)(X,Y)  \init(P,{-}\vec w)(X,Y),
  \end{array} 
 $$
 so
 $$
\begin{array}{l}
 \init(f,\vec w)(X,Y)= g(X,Y)h(X,Y)  \init(P,\vec w)(X,Y), \phantom{ssssssssssssssssss}\\
 \init(f,{-}\vec w)(X,Y) =g(X,Y)h(X,Y)  \init(P,{-}\vec w)(X,Y).
  \end{array} 
 $$
   By  condition (2) of our assumptions we see that
  $\{\,g(X,Y)=0\,\}\subset\{\,XY=0\,\}$ and
  $\{\,h(X,Y)=0\,\}\subset\{\,XY=0\,\}$.  Let us recall that the polynomials $g$ and $h$ are coprime. Using Hilbert Nullstellensatz
  we conclude that the polynomials $g$ and $h$ are powers (up to a constant) of
  different variables.  Therefore the point $(0,0)$ is the solution of the
  system $g(X,Y)=h(X,Y)=0$, which is  a~contradiction.\qed

\vspace{0.5cm} \flushright
\begin{minipage}{2.5in}
\small
Mateusz Masternak\\
Institute of Mathematics \\
Faculty of Mathematics and Natural Sciences\\
Jan Kochanowski University in Kielce \\
ul. \'Swi\e{}tokrzyska 15A\\
PL 25-406 Kielce\\
\emph{e-mail:} mateusz.masternak@ujk.edu.pl\\
\end{minipage}

\begin{thebibliography}{11}

\bibitem{A} M.F. Atiyah, {\em Angular momentum, convex polyhedra and
            algebraic geometry}, Proceedings of the Edinburgh Mathematical
            Society 26(1983) 121-138.
\bibitem{AJ} L.A. Ajzenberg, A.P. Ju\~{z}akow, {\em Integral representations
            and residues in multidimentional complex analysis\/}, (1991),
            Izdat.Nauka, Sibirskoe Otdelenie, Novosibirsk (in Russian).
\bibitem{AR} S. Abhyankar, L.A. Rubel, {\em Every difference polynomial has
           a connected zero-set\/}, J.~Indian Math. Soc., 43(1979) 69-78.
\bibitem{B}  D.N.\ Bernstein, {\em The number of roots of a system
             of equations\/},
            Funk. Anal. Appl. 9(3) (1975)),1-4 (in Russian).
\bibitem{BKuKh} D.N. Bernstein, A.G. Kouchnirenko, A.G. Khovanski, {\em
            Newton polyhedra}, UMN XXXI 3 (1976), 201-202 (in Russian).
\bibitem{P3}   Pi. Cassou-Nogu\`es, A. P\l{}oski, {\em Invariants of plane curve singularities and Newton diagrams\/},
            Univ. Iag. Acta Math Fasc. XLIX, pp. 9-34, 2011.
\bibitem{F}     W.  Fulton, {\em Algebraic curves\/}, W.A.Benjamin, Inc., (1969).
\bibitem{Kh1}   A.G. Khovanski, {\em Newton plyhedra and toroidal varieties
            \/}, Funk. Anal. Appl.,
            11(4) (1977), 56-64.
\bibitem{Kh2}   A.G. Khovanski, {\em Nevton Polyhedra and the genus of complete
            intersections\/}, Functional Anal. Appl. 12 (1978) 38-46.
\bibitem{Ku1}  A.G. Kouchnirenko, {\em Newton polyhedra and a number of roots
            of a system of $k$ equations in~$k$~variables\/}, UMN XXX
            2 (1975), 266-267 (in Russian).
\bibitem{Ku3}   A.G. Kouchnirenko, {\em Newton polyhedra and B\'ezout
              theorem\/}, Funk. Anal. Appl.,
             10(3) (1976), 82-83 (in Russian).

\bibitem{Ku4}   A.G. Kouchnirenko, {\em Poly\`{e}dres de Newton et nombres de
            Milnor\/}, Inventiones Mathematicae, 32 (1976), 1--31.
\bibitem{Ma2} M. Masternak, {\em Invariants of singularities of polynomials in two complex variables}, Universitatis Iagellonicae, Acta Mathematica, Acta Scientiarum Litterarumque MCCLV, 2001, 179-188, Krak\o{}w;
\bibitem{P1}    A. P\l{}oski, {\em Newton polygons and the \L{}ojasiewicz
            exponent of a holomorphic mapping of $\Cn^2$\/},
            Ann.\ Pol.\ Math.\ 51 (1990), 275--281.
\bibitem{P2}    A. P\l{}oski, {\em On the Irreducibility of Polynomials in Several
            Complex Variables\/}, Bulletin of the Polish Academy of Sciences Mathematics,
            Vol. 39, No. 3-4, (1991).
\bibitem{RST}   L.A.Rubel, A. Schinzel, H. Tverberg, {\em On difference polynomials and
            hereditarily irreducible polynomials\/},
            J. Number Theory, 2 (1980) 230-235.
\bibitem{We}    R. Webster, {\em Convexity\/}, Oxford University Press (1994).
\end{thebibliography}
\end{document}